\documentclass{article}
\usepackage[utf8]{inputenc}
\usepackage{amscd,epsfig,longtable}
\usepackage[T1]{fontenc}
\usepackage{verbatim}
\usepackage{amsthm}
\usepackage{amsmath}
\usepackage{amstext}
\usepackage{amssymb}
\usepackage{mathrsfs}
\usepackage{graphicx,multirow,array}
\usepackage{wasysym}
\usepackage{esint}
\usepackage{hyperref}
\usepackage[sort&compress,square,numbers]{natbib}
\usepackage{indentfirst}
\usepackage{euscript}
\usepackage{chngpage}
\usepackage{multirow}
\usepackage{enumitem}
\usepackage{xcolor}
\usepackage{subfig}
\usepackage{float}
\usepackage[labelfont=bf]{caption}
\usepackage{makecell}
\usepackage{tikz}

\RequirePackage{tikz}
\makeatletter
\def\tikzoverlay{%
    \pgfutil@ifnextchar[{\tikzoverlay@opt}{\tikzoveraly@opt[]}%
}
\def\tikzoverlay@opt[#1]#2{%
    \begin{tikzpicture}
        \node[anchor=south west, inner sep=0] (image) at (0,0) {\includegraphics[#1]{#2}};
        \newdimen\nmd@tikzoverlaywidth
        \pgfextractx{\nmd@tikzoverlaywidth}{\pgfpointanchor{image}{south east}}
        \begin{scope}
          \tikzset{x=0.01\nmd@tikzoverlaywidth}
          \tikzset{y=0.01\nmd@tikzoverlaywidth}
}
\def\endtikzoverlay{%
    \end{scope}
    \end{tikzpicture}
}
\def\tikzoverlayabs{%
    \pgfutil@ifnextchar[{\tikzoverlayabs@opt}{\tikzoveralyabs@opt[]}%
}
\def\tikzoverlayabs@opt[#1]#2{%
    \begin{tikzpicture}
        \node[anchor=south west, inner sep=0] (image) at (0,0) {\includegraphics[#1]{#2}};
        \begin{scope}
          \tikzset{x=(image.south east)}
          \tikzset{y=(image.north west)}
}
\def\endtikzoverlayabs{%
    \end{scope}
    \end{tikzpicture}
}
\makeatother

\newtheorem{theorem}{Theorem}

\newtheorem{lemma}[theorem]{Lemma}

\newtheorem{claim}{Claim}

\title{Essential Surfaces in the Exterior of $K13n586$}
\author{Chaeryn Lee}
\date{}
\begin{document}

\maketitle

\begin{abstract}
We count the number of isotopy classes of closed, connected, orientable, essential surfaces embedded in the exterior $B$ of the knot K13n586. The main result is that the count of surfaces by genus is equal to the Euler totent function. This is the first manifold for which we know the number of surfaces for any genus. The main argument is to show when normal surfaces in $B$ are connected by counting their number of components. We implement tools from Agol, Hass and Thurston to convert the problem of counting components of surfaces into counting the number of orbits in a set of integers under a collection of bijections defined on its subsets.
\end{abstract}

\section{Introduction}

Listing and counting the number of essential surfaces is one interesting approach to the study of essential surfaces in 3-manifolds. It has been shown in ~\cite{DGR} that the number of closed, orientable, essential surfaces in a given 3-manifold counted by Euler characteristic have a particular quasi-polynomial structure. For a 3-manifold $M$, let $b_M(n)$ denote the number of isotopy classes of closed, orientable, essential surfaces with Euler characteristic $n$. If $M$ is compact, orientable, irreducible, $\partial$-irreducible, atoroidal, acylindrical and does not contain any closed, nonorientable, essential surfaces then the corresponding generating function $B_M(x) = \sum_{n=1}^{\infty} b_M(-2n) x^n$ is short (Theorem 1.3 of ~\cite{DGR}).

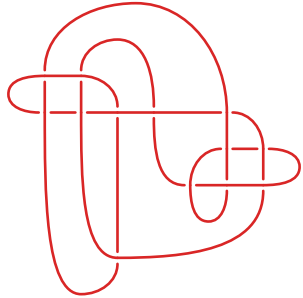
\begin{figure}[htb]
\centering
\resizebox{40mm}{!}{\definecolor{linkcolor0}{rgb}{0.85, 0.15, 0.15}
\begin{tikzpicture}[line width=2.5, line cap=round, line join=round]
  \begin{scope}[color=linkcolor0]
    \draw (7.40, 3.57) .. controls (7.40, 3.07) and (7.22, 2.55) .. 
          (6.79, 2.55) .. controls (6.29, 2.55) and (6.18, 3.18) .. (6.18, 3.76);
    \draw (6.18, 3.76) .. controls (6.18, 4.40) and (6.60, 4.97) .. (7.20, 4.97);
    \draw (7.59, 4.97) .. controls (7.86, 4.97) and (8.14, 4.97) .. (8.41, 4.97);
    \draw (8.80, 4.97) .. controls (9.30, 4.97) and (9.82, 4.80) .. 
          (9.82, 4.37) .. controls (9.82, 3.87) and (9.19, 3.76) .. (8.61, 3.76);
    \draw (8.61, 3.76) .. controls (8.20, 3.76) and (7.80, 3.76) .. (7.40, 3.76);
    \draw (7.40, 3.76) .. controls (7.06, 3.76) and (6.72, 3.76) .. (6.38, 3.76);
    \draw (5.99, 3.76) .. controls (5.08, 3.76) and (4.97, 4.94) .. (4.97, 5.99);
    \draw (4.97, 6.38) .. controls (4.97, 7.46) and (4.69, 8.61) .. 
          (3.76, 8.61) .. controls (3.13, 8.61) and (2.55, 8.19) .. (2.55, 7.59);
    \draw (2.55, 7.20) .. controls (2.55, 6.86) and (2.55, 6.52) .. (2.55, 6.18);
    \draw (2.55, 6.18) .. controls (2.55, 4.04) and (2.55, 1.34) .. (3.76, 1.34);
    \draw (3.76, 1.34) .. controls (6.06, 1.34) and (8.61, 1.60) .. (8.61, 3.57);
    \draw (8.61, 3.96) .. controls (8.61, 4.30) and (8.61, 4.64) .. (8.61, 4.97);
    \draw (8.61, 4.97) .. controls (8.61, 5.61) and (8.19, 6.18) .. (7.59, 6.18);
    \draw (7.20, 6.18) .. controls (6.46, 6.18) and (5.72, 6.18) .. (4.97, 6.18);
    \draw (4.97, 6.18) .. controls (4.57, 6.18) and (4.17, 6.18) .. (3.76, 6.18);
    \draw (3.76, 6.18) .. controls (3.43, 6.18) and (3.09, 6.18) .. (2.75, 6.18);
    \draw (2.36, 6.18) .. controls (2.08, 6.18) and (1.81, 6.18) .. (1.54, 6.18);
    \draw (1.14, 6.18) .. controls (0.65, 6.18) and (0.13, 6.36) .. 
          (0.13, 6.79) .. controls (0.13, 7.29) and (0.76, 7.40) .. (1.34, 7.40);
    \draw (1.34, 7.40) .. controls (1.75, 7.40) and (2.15, 7.40) .. (2.55, 7.40);
    \draw (2.55, 7.40) .. controls (3.19, 7.40) and (3.76, 6.98) .. (3.76, 6.38);
    \draw (3.76, 5.99) .. controls (3.76, 4.50) and (3.76, 3.02) .. (3.76, 1.54);
    \draw (3.76, 1.14) .. controls (3.76, 0.55) and (3.19, 0.13) .. 
          (2.55, 0.13) .. controls (1.34, 0.13) and (1.34, 3.56) .. (1.34, 6.18);
    \draw (1.34, 6.18) .. controls (1.34, 6.52) and (1.34, 6.86) .. (1.34, 7.20);
    \draw (1.34, 7.59) .. controls (1.34, 8.99) and (2.83, 9.82) .. 
          (4.37, 9.82) .. controls (6.16, 9.82) and (7.40, 8.09) .. (7.40, 6.18);
    \draw (7.40, 6.18) .. controls (7.40, 5.78) and (7.40, 5.38) .. (7.40, 4.97);
    \draw (7.40, 4.97) .. controls (7.40, 4.64) and (7.40, 4.30) .. (7.40, 3.96);
  \end{scope}
\end{tikzpicture}}
\caption{The knot $K13n586$.}
\end{figure}

No general results are known when we turn to counting \textit{connected} essential surfaces. Define $a_M(g)$ as the number of isotopy classes of closed, connected, orientable surfaces of genus $g$ in $M$. To relate this index to Euler characteristic $n = g - 1$, we may adjust it and define $\tilde{a}_M(n) = a_M(n + 1)$. Section 8 of ~\cite{DGR} provides data of $a_M(g)$ for $g \le 200$ of certain manifolds but there are very few manifolds for which we actually know the full list of essential surfaces embedded inside. The only known cases are where $a_M(g) = 0$ for all large $g$. In this paper we study a specific manifold, the exterior of the knot $K13n586$, and present a complete description of all of its essential connected surfaces counted by Euler characteristic. This is the first manifold for which we know $a_M(g)$ for all $g$. In particular, we will prove the following result which was asserted as a conjecture in Section 8.1 of ~\cite{DGR}. 

\begin{theorem}
\label{main theorem}
For the exterior $B$ of the knot $K13n586$, $\tilde{a}_B(n) = \phi(n)$, $n > 1$ where $\phi(n)$ is the Euler totent function.
\end{theorem}

In Section 2, we give a detailed description of the triangulation of $B$ and the normal surfaces embedded in $B$ with respect to this triangulation. In Section 3, we look at combinatorial tools that enable us to count the number of components of surfaces. We then prove our theorem in Section 4.

This work was partially supported by U.S. National Science Foundation grant DMS-1811156.

\section{The structure of $B$ and its normal surfaces}

We first look at some normal surface theory introduced in the work by ~\cite{T}. Let $M$ be a manifold with a triangulation $\mathcal{T}$. A surface is said to be \textit{normal} if it is in general position with the 1-skeleton $\mathcal{T}^{(1)}$ of $\mathcal{T}$ and meets each tetrahedron only in triangles and quadrilaterals that we call elementary disks. An \textit{elementary disk} $E$ in a tetrahedron $\Delta$ of $\mathcal{T}$ is a properly embedded disk that meets each edge of $\Delta$ in at most one point and each face of $\Delta$ in at most one line. We will follow the convention that if $E \cap \mathcal{T}^{(1)}$ is a planar set then $E$ is planar and if not then $E$ is the cone $b \ast \partial E$ where $b$ is the centroid of the 3-simplex $\Delta$ spanned by $E \cap \mathcal{T}^{(1)}$. Elementary disks and hence normal surfaces are uniquely determined by their intersection points with $\mathcal{T}^{(1)}$. A \textit{normal isotopy} of $M$ is an isotopy that leaves every simplex of $\mathcal{T}$ invariant. There are exactly 7 normal isotopy classes of elementary disks, four triangles that each cut off a vertex and three quadrilaterals that separate the three pairs of disjoint edges. We call such normal isotopy classes the \textit{disk types} of a tetrahedron. Suppose we have two normal surfaces $S$ and $T$ that intersect transversely and have the same quad types in each tetrahedron of $\mathcal{T}$. Then there is a  unique normal surface determined by $(S \cup T) \cap \mathcal{T}^{(1)}$ which we call the \textit{normal sum} of $S$ and $T$ and write $S + T$.

\begin{figure}[h]
\begin{center}
\begin{tikzoverlay}[width=74mm]{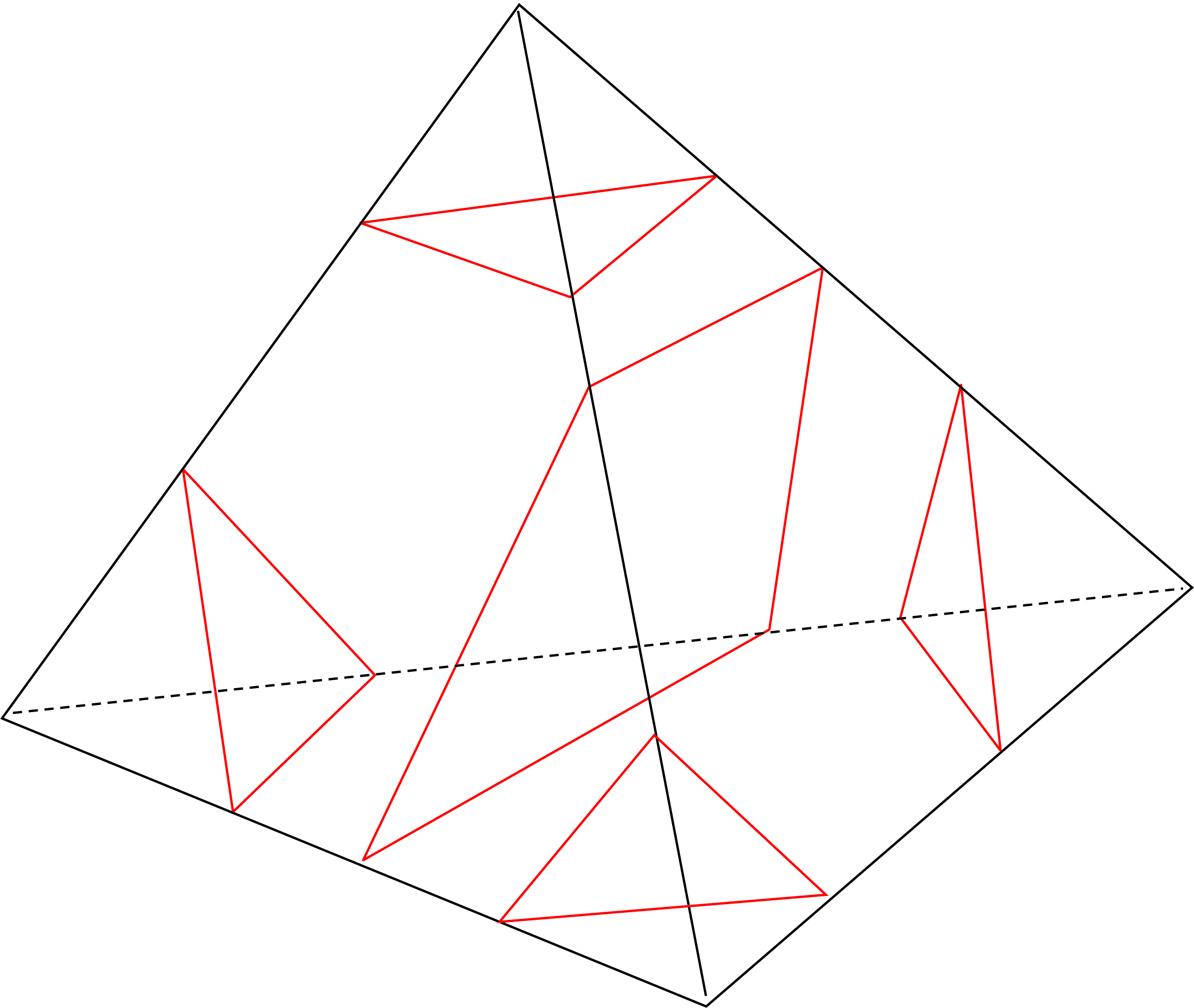}
\tikzstyle{every node}=[font=\small]
    \node[above] at (43.2,85) {$0$};
    \node[left] at (-0.5,24) {$1$};
    \node[below] at (59.5,-0.5) {$2$};
    \node[right] at (100.5,35) {$3$};
    
    \node[left] at (22,55) {E01};
    \node[left] at (53,37) {E02};
    \node[right] at (73,60) {E03};
    \node[left] at (26,12) {E12};
    \node[below] at (70,32) {E13};
    \node[below] at (78.5,15) {E23};
    
    \node[above,red] at (41,68) {tri0};
    \node[left,red] at (17,34) {tri1};
    \node[below,red] at (54,9) {tri2};
    \node[right,red] at (85,40) {tri3};
    \node[red] at (60,50) {quad};
    
    \draw [-stealth](50,-10) -- (50,-5);
    \node[below] at (50,-10) {F0};
    \draw [-stealth](87,70) -- (83,67);
    \node[above right] at (87,70) {F1};
    \draw [-stealth](50,95) -- (50,90);
    \node[above] at (50,95) {F2};
    \draw [-stealth](6,63) -- (10,60);
    \node[above left] at (6,63) {F3};
\end{tikzoverlay}
\end{center}
  \caption{Standard tetrahedron with elementary disk types, labelled with names of vertices, edges and faces.}
  \label{fig:fig0}
\end{figure}
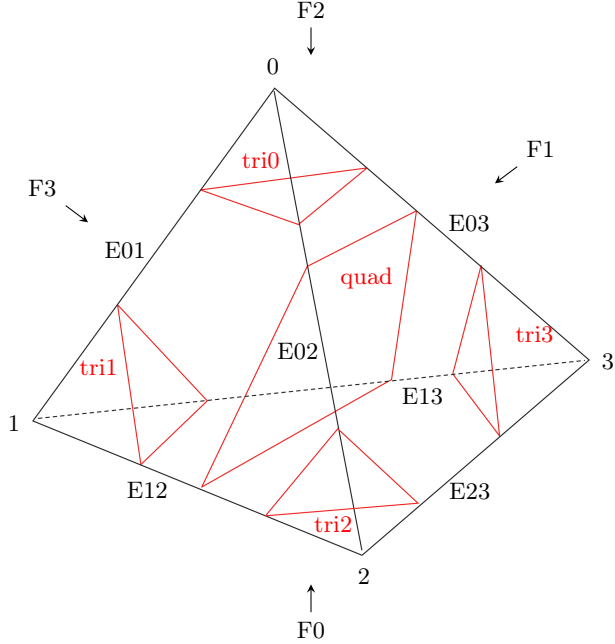

\begin{table}[h]
\centering
{\small
\begin{tabular}{c|cccc}
\noalign{\smallskip}\noalign{\smallskip}\hline\hline
& F0 & F1 & F2 & F3 \\
\hline
tet0 & tet2 & tet1 & tet7 & tet5 \\
 & (0, 1, 3, 2) &(0, 1, 3, 2) &(0, 2, 1, 3) &(2, 1, 0, 3) \\
tet1 & tet4 & tet0 & tet6 & tet3 \\ 
 & (3, 0, 1, 2) & (0, 1, 3, 2) & (3, 1, 2, 0) & (0, 1, 3, 2) \\
tet2 & tet0 & tet5 & tet3 & tet9  \\
& (0, 1, 3, 2) & (0, 2, 1, 3) & (1, 3, 0, 2) & (3, 1, 2, 0)  \\
tet3 & tet2 & tet9 & tet1 & tet8  \\
& (2, 0, 3, 1) & (0, 1, 3, 2) & (0, 1, 3, 2) & (2, 1, 0, 3)  \\
tet4 & tet7 & tet8 & tet9 & tet1  \\
& (3, 1, 2, 0) & (0, 1, 3, 2) & (3, 1, 2, 0) & (1, 2, 3, 0)  \\
tet5 & tet7 & tet6 & tet2 & tet0  \\
& (0, 1, 3, 2) & (0, 1, 3, 2) & (0, 2, 1, 3) & (2, 1, 0, 3)  \\
tet6 & tet9 & tet5 & tet1 & tet8  \\
& (3, 0, 1, 2) & (0, 1, 3, 2) & (3, 1, 2, 0) & (0, 1, 3, 2)  \\
tet7 & tet5 & tet0 & tet8 & tet4  \\
& (0, 1, 3, 2) & (0, 2, 1, 3) & (1, 3, 0, 2) & (3, 1, 2, 0)  \\
tet8 & tet7 & tet4 & tet6 & tet3  \\
& (2, 0, 3, 1) & (0, 1, 3, 2) & (0, 1, 3, 2) & (2, 1, 0, 3)  \\
tet9 & tet2 & tet3 & tet4 & tet6  \\
& (3, 1, 2, 0) & (0, 1, 3, 2) & (3, 1, 2, 0) & (1, 2, 3, 0)  \\
\hline
\hline
\end{tabular}
\caption {Triangulation of $\mathcal{S}$ of $B$. The $(i, j)$-th entry shows which tetrahedron is glued to face F$j - 1$ of tetrahedron tet$i - 1$. The permutations below show the order with which their vertices are glued.}
\label{table:tb1}
}
\end{table}

\begingroup
\renewcommand{\arraystretch}{1.2}
\begin{table}
\centering
{\small
\begin{tabular}{c|c}
\noalign{\smallskip}\noalign{\smallskip}\hline\hline
& Edges Identified \\
\Xhline{2\arrayrulewidth}
e0 & E10 : tet0, E23 : tet4, E12 : tet5, E20 : tet7, E13 : tet7, E32: tet8 \\
\Xhline{0.0001\arrayrulewidth}
e1 & E20 : tet0, E30 : tet1, E02 : tet5, E03 : tet6 \\
\Xhline{0.0001\arrayrulewidth}
e2 & E12 : tet0, E20 : tet2, E13 : tet2, E32 : tet3, E10 : tet5, E23 : tet9 \\
\Xhline{0.0001\arrayrulewidth}
\multirow{2}{*}{e3} & E30 : tet0, E01 : tet1, E20 : tet1, E01 : tet3, E30 : tet3 \\
& E31 : tet6, E30 : tet7, E21 : tet8, E20 : tet9 \\
\Xhline{0.0001\arrayrulewidth}
e4 & E31 : tet0, E21 : tet2, E23 : tet5, E32 : tet6, E32 : tet7, E21 : tet9 \\
\Xhline{0.0001\arrayrulewidth}
e5 & E23 : tet0, E32 : tet1, E32 : tet2, E21 : tet4, E31 : tet5, E21 : tet7 \\
\Xhline{0.0001\arrayrulewidth}
e6 & E12 : tet1, E01 : tet2, E13 : tet3, E01 : tet4, E31 : tet9 \\
\Xhline{0.0001\arrayrulewidth}
\multirow{2}{*}{e7} & E31 : tet1, E03 : tet2, E21 : tet3, E20 : tet4, E30 : tet5 \\
&  E01 : tet6, E20 : tet6, E01 : tet8, E30 : tet8 \\
\Xhline{0.0001\arrayrulewidth}
e8 & E20 : tet3, E03 : tet4, E02 : tet8, E30 : tet9 \\
\Xhline{0.0001\arrayrulewidth}
e9 & E31 : tet4, E12 : tet6, E01 : tet7, E13 : tet8, E01 : tet9 \\ 
\hline
\hline
\end{tabular}
\caption {List of edges in each tetrahedron that is glued to the edge e$i$ of $\mathcal{S}$.}
\label{table:tb2}
}
\end{table}
\endgroup

\begin{table}[h]
\centering
{\small
\begin{tabular}{c|ccccc|ccccc}
\noalign{\smallskip}\noalign{\smallskip}\hline\hline
\multirow{2}{*}{} & \multicolumn{5}{c|}{$F$} & \multicolumn{5}{c}{$G$} \\
\cline{2-11}
      & tri0  & tri1 & tri2 & tri3 & quad & tri0 & tri1 & tri2 & tri3 & quad \\
\hline
 tet0 & 2 & 0 & 0 & 0 & 2 & 2 & 0 & 2 & 0 & 0 \\
 tet1 & 3 & 1 & 0 & 0 & 1 & 1 & 1 & 0 & 2 & 1 \\
 tet2 & 1 & 1 & 0 & 0 & 1 & 2 & 0 & 0 & 2 & 0 \\
 tet3 & 3 & 1 & 1 & 1 & 0 & 1 & 1 & 2 & 0 & 1 \\
 tet4 & 2 & 0 & 0 & 2 & 0 & 2 & 0 & 1 & 1 & 1 \\
 tet5 & 2 & 0 & 2 & 0 & 0 & 2 & 0 & 0 & 0 & 2 \\
 tet6 & 1 & 1 & 0 & 2 & 1 & 3 & 1 & 0 & 0 & 1 \\
 tet7 & 2 & 0 & 0 & 2 & 0 & 1 & 1 & 0 & 0 & 1 \\
 tet8 & 1 & 1 & 2 & 0 & 1 & 3 & 1 & 1 & 1 & 0 \\
 tet9 & 2 & 0 & 1 & 1 & 1 & 2 & 0 & 0 & 2 & 0 \\
\hline
\hline
\end{tabular}

\begin{tabular}{c|ccccc}
\noalign{\smallskip}\noalign{\smallskip}\hline\hline
\multirow{2}{*}{} & \multicolumn{5}{c}{$uF + vG$} \\
\cline{2-6}
      & tri0  & tri1 & tri2 & tri3 & quad \\
\hline
 tet0 & $2u + 2v$ & 0 & $2v$ & 0 & $2u$ \\
 tet1 & $3u + v$ & $u + v$ & 0 & $2v$ & $u + v$ \\
 tet2 & $u + 2v$ & $u$ & 0 & $2v$ & $u$ \\
 tet3 & $3u + v$ & $u + v$ & $u + 2v$ & $u$ & $v$ \\
 tet4 & $2u + 2v$ & 0 & $v$ & $2u + v$ & $v$ \\
 tet5 & $2u + 2v$ & 0 & $2u$ & 0 & $2v$ \\
 tet6 & $u + 3v$ & $u + v$ & 0 & $2u$ & $u + v$ \\
 tet7 & $2u + v$ & $v$ & 0 & $2u$ & $v$ \\
 tet8 & $u + 3v$ & $u + v$ & $2u + v$ & $v$ & $u$ \\
 tet9 & $2u + 2v$ & 0 & $u$ & $u + 2v$ & $u$ \\
\hline
\hline
\end{tabular}
\caption{Number of each disk type that forms the surfaces $F$, $G$ and $uF + vG$. Disk types are depicted in Figure {\ref{fig:fig0}}.}
\label{table:tb3}
}
\end{table}

Throughout this paper we will fix a specific triangulation $\mathcal{S}$ of $B$ that is provided by snappy.HTLinkExteriors on SnapPy ~\cite{CDGW} (the current labelling on the tetrahedra differs). Tables~\ref{table:tb1} and ~\ref{table:tb2} give a description of this triangulation. The triangulation $\mathcal{S}$ consists of a single vertex, 10 edges, named e$i$ ($0 \le i \le 9$), and 10 tetrahedra, named tet$j$ $0 \le i \le 9$, each with an ordering $0, 1, 2, 3$ on their vertices as in Figure {\ref{fig:fig0}}. In each tetrahedron E$ij$ denotes the edge between vertices $i$ and $j$ while F$k$ denotes the face opposite to vertex $k$, that is the 2-simplex spanned by all vertices but $k$. Note that an orientation is fixed on edges E$ij$ from $i$ to $j$ so that E$ij$ and E$ji$ represent the same edge but with opposite orientation. The entry in Table~\ref{table:tb1} in the row labeled tet$i$ and column labeled F$j$ shows which tetrahedron is glued to face F$(j - 1)$ of tetrahedron tet$(i - 1)$ (the discrepancy in the indices comes from the fact that labels in our triangulation start at 0) and the order with which their vertices are glued. For instance, the vertices $(0, 1, 3)$ in F$2$ of tet$2$ are glued in order to the vertices $(1, 3, 2)$ which is F$0$ of tet$3$. Note that we give the full permutations of size 4 in Table~\ref{table:tb1} but one value in each permutation is redundant since every face is made up of 3 vertices. Table~\ref{table:tb2} gives a list of the edges of the 10 tetrahedra that are identified to each $e_i$. The orientation with which they are identified are taken into account by the indices E$ij$.

Section 8.1 of ~\cite{DGR} shows that there is a bijection from the isotopy classes of essential surfaces in $B$ to surfaces of the form $uF + vG$ ($u, v \in \mathbb{N}$). In our triangulation $\mathcal{S}$ the only disk types we consider are the triangles, ordered by the vertex they are positioned at, and the quadrilateral that separates the edges E01 and E23 as in Figure {\ref{fig:fig0}}. Note that $F$ and $G$ have the same disk types, hence it is possible to define their normal sum. Table~\ref{table:tb3} gives a description of the normal surfaces $F$, $G$ and how they are embedded in $B$ with respect to $\mathcal{S}$. Table~\ref{table:tb3} gives the number of each disc type that forms the surfaces $F$, $G$ and $uF + vG$. Both $F$ and $G$ are surfaces of genus 2.

Our goal is to count $\tilde{a}_B(n)$, the number of connected surfaces $uF + vG$ with genus $n + 1$ or Euler characteristic $-2n$. Following the work of ~\cite{DGR} we can in fact reduce Theorem~\ref{main theorem} to Lemma~\ref{main lemma} below. Since $F$ and $G$ are surfaces of genus 2, surfaces $uF + vG$ with Euler characteristic $-2n$ correspond to lattice points in $\mathbb{N}^2$ on the line $x + y = n$. We want to show that $\tilde{a}_B(n)$ is exactly the number of all primitive such lattice points. If $uF + vG$ is connected then $gcd(u,  v) = 1$, otherwise $uF + vG$ would consist of disjoint copies of the surface $\frac{1}{gcd(u,  v)}(uF + vG)$. Hence it suffices to show the converse to prove Theorem~\ref{main theorem}.

\begin{lemma}
\label{main lemma}
If $gcd(u,  v) = 1$ then $uF + vG$ is a connected surface.
\end{lemma}

The rest of this paper will focus on proving  Lemma~\ref{main lemma}.

\section{Counting orbits of integer intervals}
In this section we follow the work of Section 4 in ~\cite{AHT} and introduce some techniques used to count the number of orbits in a set of integers under a particular pseudogroup action. Fixing an order and orientation on the edges of our triangulation $\mathcal{S}$, we may correspond an integer interval $[1, N] \subseteq \mathbb{N}$ to the intersection points of $uF + vG$ with the edges of $\mathcal{S}$. The intersection arcs of $uF + vG$ with the faces of $\mathcal{S}$ then give rise to bijections defined on subintervals of $[1, N]$. Hence we can change the problem of counting the number of components of $uF + vG$ to counting the number of orbits in the interval $[1, N]$ under these bijections.

Let  $[1, N] \subseteq \mathbb{N}$ be a set of integers. A bijection $g\colon [a, b] \to [c, d]$ defined on subintervals  $[a, b], [c, d] \subseteq [1, N]$ is called a \textit{pairing}. A pairing is \textit{orientation preserving} if it is increasing and \textit{orientation reversing} if it is decreasing. The \textit{width} of a pairing $g\colon [a, b] \to [c, d]$ is the number of integers in its domain or range, and denoted $|g| = b - a + 1 = d - c + 1$. Note that all sets and maps we work with henceforth are discrete. 

Assume that we have the interval $[1, N]$ and a collection of pairings $\{ g_i \}, 1 \le i \le k$ defined on subintervals of $[1, N]$. These pairings naturally generate a pseudogroup on $[1, N]$. Any two integers are said to be in the same \textit{orbit} if some pairing in the pseudogroup generated by $\{ g_i \}$ sends one to the other. Note that orbits form a partition of $[1, N]$. Our goal is to count the number of such orbits. 

We now describe some modifications we can perform on the interval $[1, N]$ and pairings $\{ g_i \}$. An important trait of these modifications is that they do not change the number of orbits (Lemma 11 of ~\cite{AHT}).

\begin{description}[style=nextline]
 \item[Trimming]
 Trimming modifies an orientation reversing pairing so that its domain and range are disjoint. Let $g\colon [a, b] \to [c, d]$ be an orientation reversing pairing so that its domain and range overlap, $g(a) = d, g(b) = c, c \le b$. Define the pairing $g'\colon [a, \frac{a + d}{2}) \to (\frac{a + d}{2}, d]$ and replace $g$ with $g'$. We say that $g'$ is obtained by \textit{trimming} $g$. $g'$ is also an orientation reversing pairing.    
 
 \item[Truncation]
 Suppose there is some rightmost subinterval $[N' + 1, N]$ of $[1, N]$ that lies in the range of only one pairing $g\colon [a, b] \to [c, N]$, $c \le N' + 1 \le N$. Truncation "peels off" this subinterval from $g$ and shortens the entire interval $[1, N]$ to  $[1, N']$ without changing the orbit structure. If $g$ is orientation preserving define $g'\colon [a, b - (N - N')] \to [c, N']$. If $g$ is orientation reversing first make sure that its domain and range are disjoint, trim $g$ if necessary. Then define $g'\colon [a + (N - N'), b] \to [c, N']$. Replacing $g$ with $g'$  is called \textit{truncating} $g$.
 
 We can similarly truncate a leftmost subinterval $[1, N'']$ of $[1, N]$ that lies in the domain and range of a single pairing. After this modification is applied the interval $[1, N]$ is shortened to $[N'' + 1, N]$. To ensure our interval starts at 1 we may translate the domain and range of every pairing in $\{ g_i \}$ by $-N''$. 

 \item[Transmission]
 Transmission composes two pairings to shift the domain and range of one by the other. Assume $g_1$ has disjoint domain and range and that the range of $g_2$ lies in the range of $g_1$. If the domain of $g_2$ does not lie in the range of $g_1$, define $g_2' = {g_1}^{-1} \circ g_2$. If the domain of $g_2$ does lie in the range of $g_1$, define $g_2' = {g_1}^{-1} \circ g_2 \circ g_1$. Replacing $g_2$ with $g_2'$ we \textit{transmit} $g_2$ \textit{by} $g_1$.        
 
\end{description}

\section{Proof of Lemma~\ref{main lemma}}
We will order the edges of $\mathcal{S}$ as $e1, e8, e6, e9, e0, e2, e3, e4, e5, e7$ (this order is intentionally chosen in this way) then label each of the intersection points by the integers $1, 2, ..., 24u + 24v$ to obtain an interval $[1, 24u + 24v] \subseteq \mathbb{N}$. Faces of $\mathcal{S}$ intersect the surface $uF + vG$ in a collection of arcs joining edges of each face. A collection of arcs that join the same edges of a face will be called an \textit{arc type}. Every arc type defines a pairing $p_i\colon [a_i, b_i] \to [c_i, d_i]$ where $[a_i, b_i]$ and $[c_i, d_i]$ are subintervals of $[1, 24u + 24v]$. To prove Lemma~\ref{main lemma} our goal is to show that there is one orbit of $[1, 24u + 24v]$ under the collection of pairings $\{p_i\}$.

\begin{claim} 
The number of orbits of interval $[1, 24u + 24v]$ under pairings $\{p_i\}$ is at most the number of orbits of the interval $[1, 2u + 2v]$ under the orientation reversing pairings 
\begin{align*}
& \{ f\colon [1, u] \to [u + 1, 2u], \\
& g\colon [2u + 1, 2u + v] \to [2u + v + 1, 2u + 2v], \\ 
& h\colon [1, u + v] \to [u + v + 1, 2u + 2v] \}.
\end{align*}

\begin{figure}[h]
\begin{center}
\begin{tikzoverlay}[width=30mm]{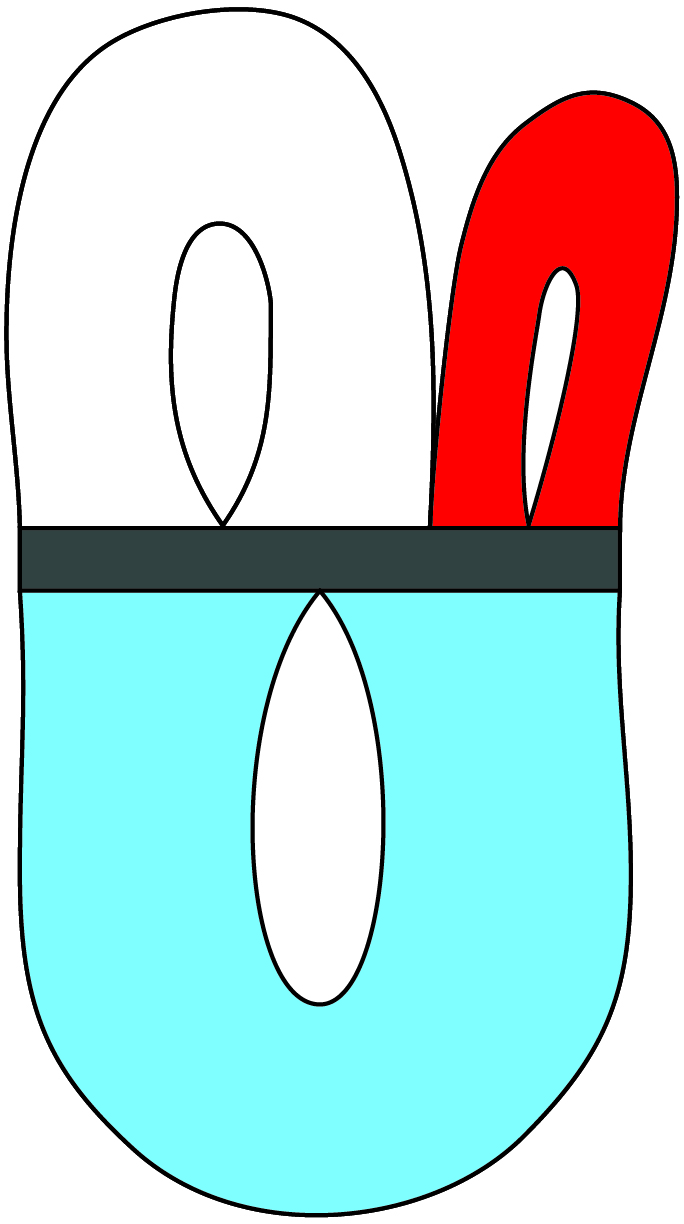}
\tikzstyle{every node}=[font=\small]

    \node[above] at (10,170) {$f$};
    \node[above,red] at (94,165) {$g$};
    \node[below,cyan] at (70,5) {$h$};

\end{tikzoverlay}
\end{center}
  \caption{The thick interval in the middle represents the interval $[1, 2u+2v]$. Pairings $f, g$ and $h$ are depicted as the bands attached to this interval.}
  \label{fig:fig4}
\end{figure}

\end{claim}

\begin{proof}

The first step will be to reduce the interval $[1, 24u + 24v]$ to its subinterval corresponding to $e1$ which is $[1, 4u + 4v]$. We start from the rightmost subinterval corresponding to $e7$. Observe that all pairings that have their ranges on the edge $e7$ (these include the pairings that have both their domain and range on $e7$) can be transmitted by the pairing $t_7$ depicted as in Figure {\ref{fig:fig1}} so that their ranges are shifted into $e1$. This is because the domain of $t_7$ lies on $e1$ and its range is the entire $e7$. For any pairing whose domain lies on $e7$ and range lies on a distinct edge, simply swap the domain and range and apply transmission. Now $e7$ only lies in the range of $t_7$ so we can truncate $t_7$ and remove the subinterval corresponding to $e7$ from $[1, 24u + 24v]$. Similarly, all pairings that have their range on the edges $e0, e2, e3, e4, e5$ can be transmitted by the pairings $t_0, t_2, t_3, t_4, t_5$ into $e_1$ and $t_0, t_2, t_3, t_4, t_5$ can be truncated. All pairings defined on $e6, e9$ can be transmitted by pairings $t_6, t_9$ into $e8$ and $t_6, t_9$ are truncated. Finally, all pairings defined on $e8$ can be transmitted by $t_{8a}, t_{8b}$ into $e1$ where pairings whose range lie in both the range of $t_{8a}, t_{8b}$ can be cut into two pairings before transmission.

\begin{figure}[h]
\begin{center}
\begin{tikzoverlay}[width=100mm]{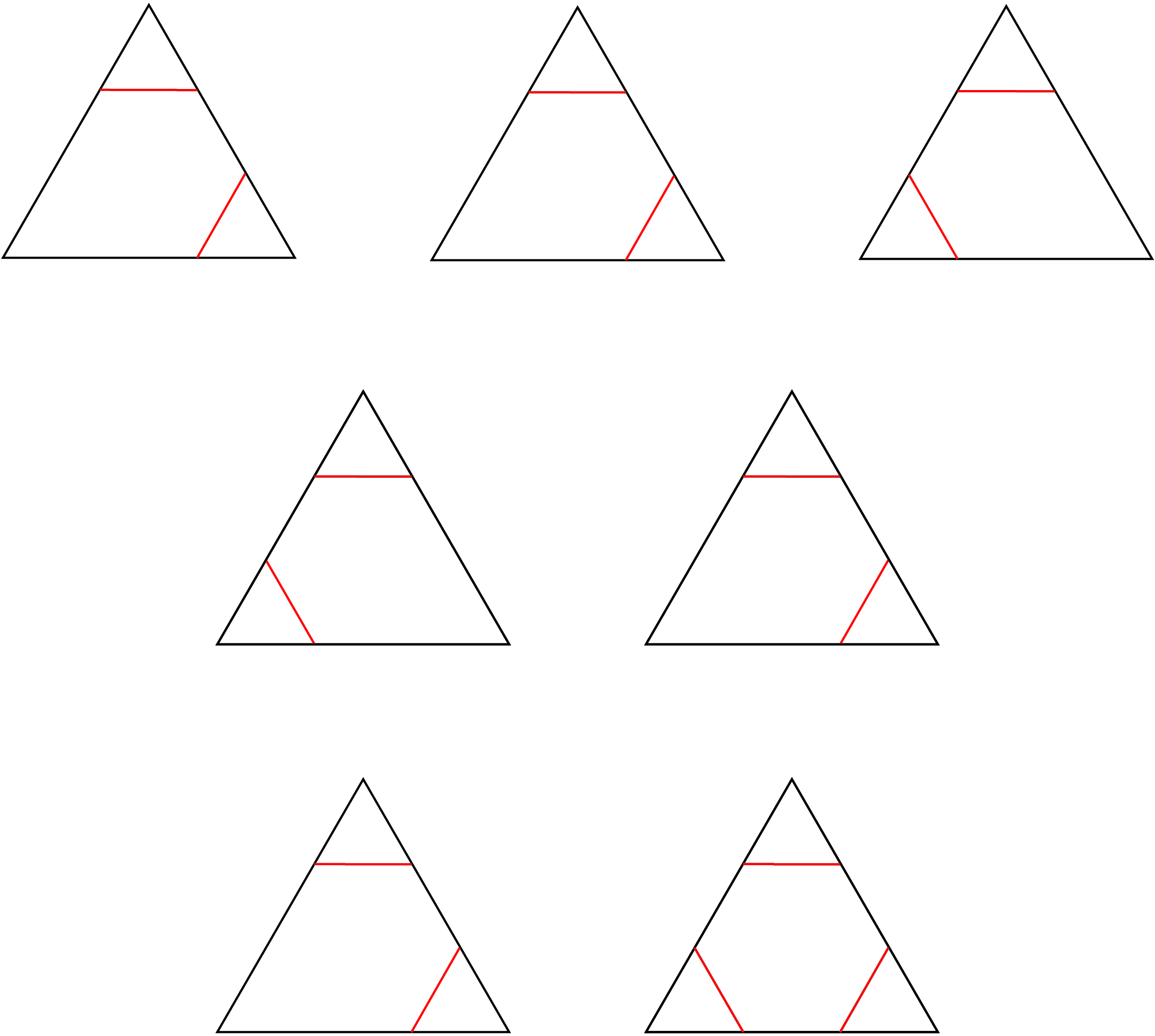}
\tikzstyle{every node}=[font=\small]
    \node[left] at (0,67) {$1$};
    \node[right] at (26,67) {$2$};
    \node[above] at (13,89.5) {$0$};
    \node[right] at (19.5,78) {e1};
    \node[below] at (13,67.5) {e0};
    \node[below] at (13,64) {F3 of tet5};
    \node[left,red,font=\footnotesize] at (20,72) {$t_0 (2u+2v)$};
    
    \node[left] at (37,67) {$1$};
    \node[right] at (63,67) {$2$};
    \node[above] at (50,89.5) {$0$};
    \node[right] at (56.5,78) {e1};
    \node[below] at (50,67.5) {e2};
    \node[below] at (50,64) {F3 of tet0};
    \node[left,red,font=\footnotesize] at (57,72) {$t_2 (2u+2v)$};
    
    \node[left] at (74,67) {$2$};
    \node[right] at (100,67) {$3$};
    \node[above] at (87,89.5) {$0$};
    \node[left] at (80.5,78) {e1};
    \node[right] at (93.5,78) {e3};
    \node[below] at (87,67.5) {e5};
    \node[below] at (87,64) {F1 of tet0};
    \node[below,red,font=\footnotesize] at (87,82) {$t_3$};
    \node[below,red,font=\footnotesize] at (87,79.5) {$(4u+2v)$};
    \node[right,red,font=\footnotesize] at (80.5,72) {$t_5 (2v)$};
 
    \begin{scope}[shift={(18.5,0)}]
    \node[left] at (0,33) {$2$};
    \node[right] at (26,33) {$3$};
    \node[above] at (13,55.5) {$0$};
    \node[left] at (6.5,44) {e1};
    \node[right] at (19.5,44) {e7};
    \node[below] at (13,34) {e4};
    \node[below] at (13,30.5) {F1 of tet5};
    \node[right,red,font=\footnotesize] at (6.5,38) {$t_4 (2u)$};
    \node[below,red,font=\footnotesize] at (13,49) {$t_7$};
    \node[below,red,font=\footnotesize] at (13,46.5) {$(2u+4v)$};
    \end{scope}

    \begin{scope}[shift={(-18.5,0)}]
    \node[left] at (74,33) {$1$};
    \node[right] at (100,33) {$3$};
    \node[above] at (87,55.5) {$0$};
    \node[right] at (93.5,44) {e8};
    \node[below] at (87,34) {e6};
    \node[below] at (87,30.5) {F2 of tet9};
    \node[left,red,font=\footnotesize] at (93.5,38) {$t_6 (2u+2v)$};
    \end{scope}
    
    \node[left] at (18.5,-0.5) {$1$};
    \node[right] at (44.5,-0.5) {$3$};
    \node[above] at (31.5,22) {$0$};
    \node[right] at (38,11) {e8};
    \node[below] at (31.5,0.5) {e9};
    \node[below] at (31.5,-3) {F2 of tet4};
    \node[left,red,font=\footnotesize] at (38,5) {$t_9 (2u+2v)$};
    
    \node[left] at (55.5,-0.5) {$1$};
    \node[right] at (81.5,-0.5) {$2$};
    \node[above] at (68.5,22) {$0$};
    \node[left] at (62,11) {e3};
    \node[right] at (75,11) {e8};
    \node[below] at (68.5,0.5) {e7};
    \node[below] at (68.5,-3) {F3 of tet3};
    \node[left,red,font=\footnotesize] at (76,6.5) {$t_{8b}$};
    \node[left,red,font=\footnotesize] at (75,3.5) {$(u+3v)$};
    \node[below,red,font=\footnotesize] at (68.5,15) {$t_{8a}$};
    \node[below,red,font=\footnotesize] at (68.5,12.5) {$(3u+v)$};
\end{tikzoverlay}
\end{center}
  \caption{Faces with vertices, edges, and pairings $t_i$ used for transmission. The numbers in parentheses are the widths of each pairing. Details about pairings other than $t_i$ are omitted.}
  \label{fig:fig1}
\end{figure}
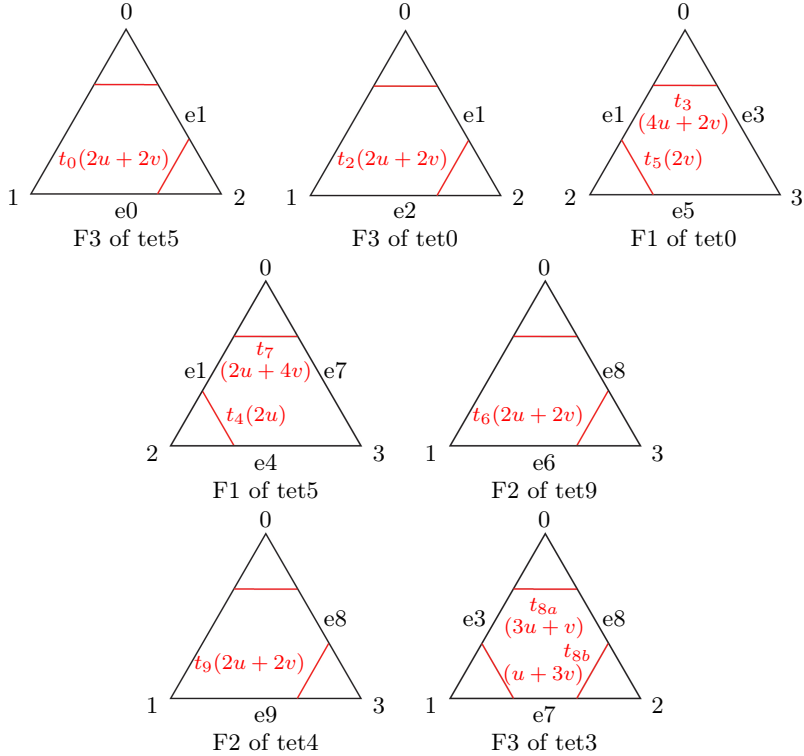

We now have $[1, 4u + 4v]$ and pairings $\{p'_i\}$ where the domains and ranges of $p'_i$ are all subintervals of $[1, 4u + 4v]$. Let $\{p'_{i1}, p'_{i2}, ..., p'_{in} \}$ be a subset of $\{p'_i\}$. If there are $n$ orbits of $[1, 4u + 4v]$ under $\{p'_{i1}, p'_{i2}, ..., p'_{in} \}$ there are at most $n$ orbits of $[1, 4u + 4v]$ under $\{p'_i\}$.  

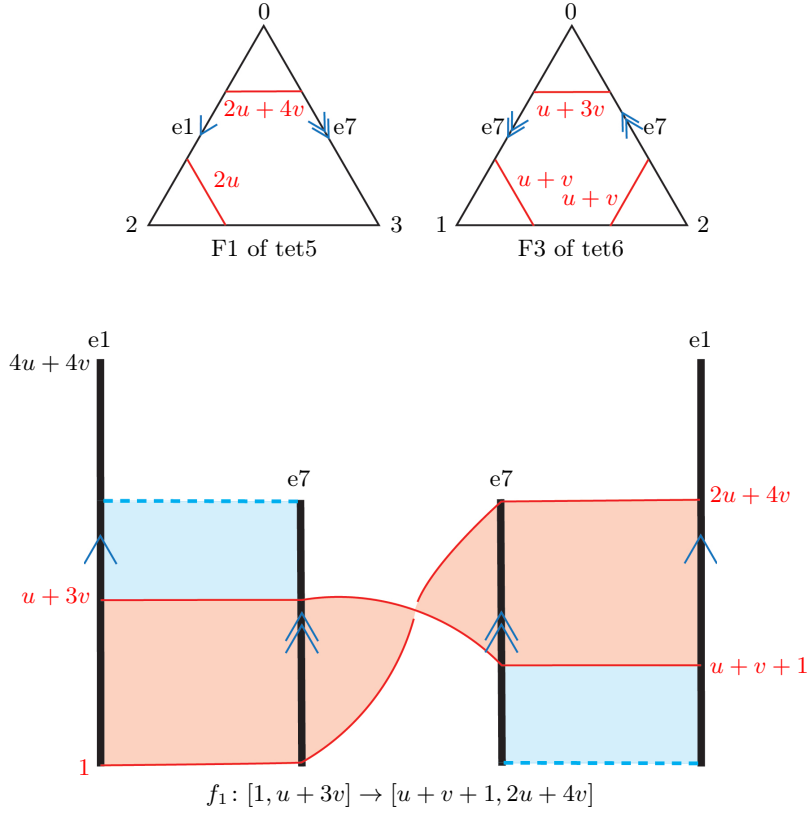
\begin{figure}[H]
\begin{center}
\begin{tikzoverlay}[height=100mm]{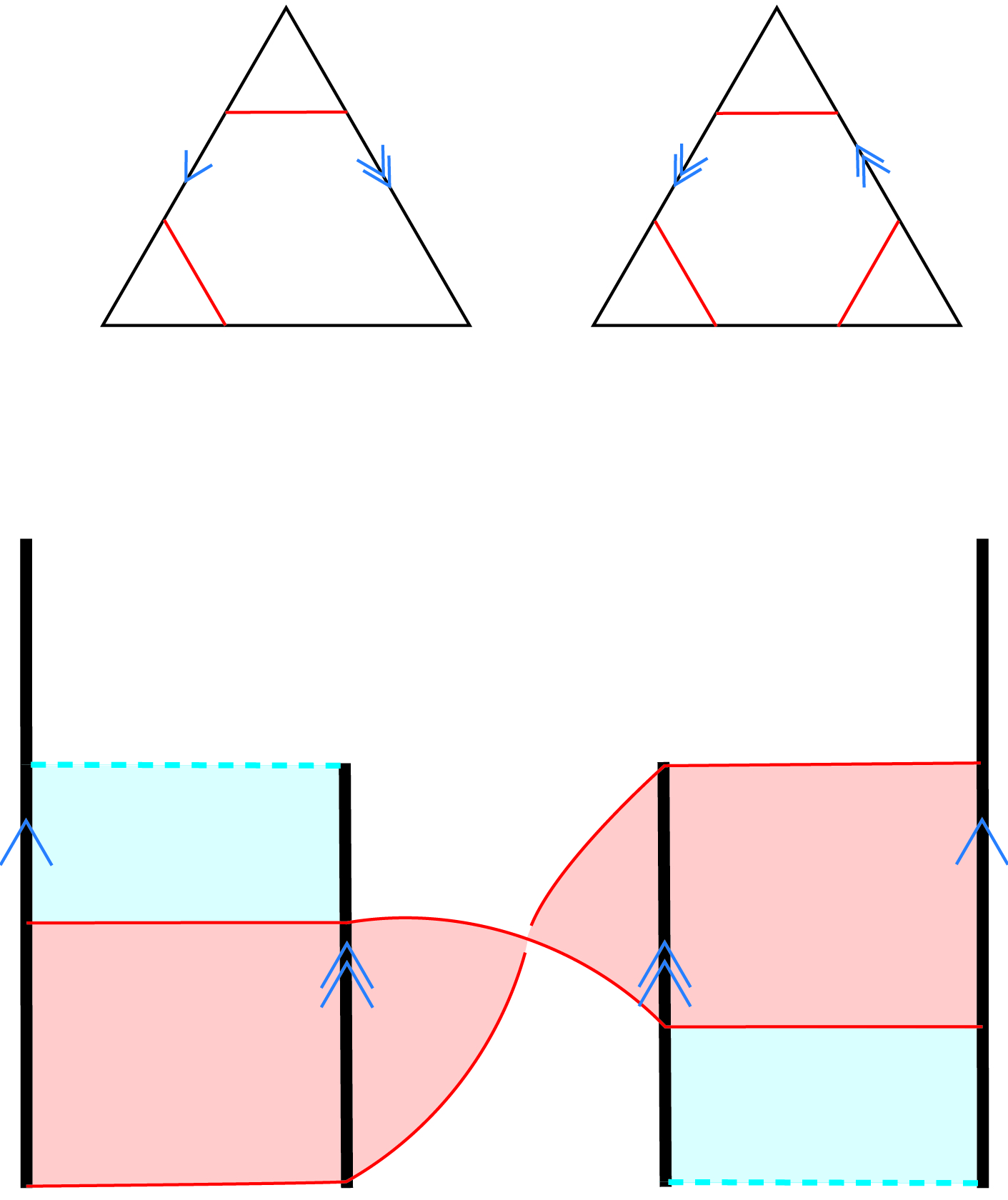}
\tikzstyle{every node}=[font=\small]
    \node[left] at (10,87) {$2$};
    \node[right] at (47,87) {$3$};
    \node[above] at (28.5,118) {$0$};
    \node[left] at (19,102.5) {e1};
    \node[right] at (38,102.5) {e7};
    \node[below] at (28.5,86) {F1 of tet5};
    \node[below,red] at (28.5,108) {$2u+4v$};
    \node[right,red] at (19,94) {$2u$};
    
    \node[left] at (59,87) {$1$};
    \node[right] at (95.5,87) {$2$};
    \node[above] at (77,118) {$0$};
    \node[left] at (68,102.5) {e7};
    \node[right] at (87,102.5) {e7};
    \node[below] at (77,86) {F3 of tet6};
    \node[below,red] at (77, 108) {$u+3v$};
    \node[right,red] at (67,94) {$u+v$};
    \node[left,red] at (86,91) {$u+v$};
    
    \node[above] at (2.5,66) {e1};
    \node[above] at (34,44) {e7};
    \node[above] at (66,44) {e7};
    \node[above] at (97.5,66) {e1};
    
    \node[left] at (2.5,65) {$4u+4v$};
    \node[left,red] at (2.5,1) {$1$};
    \node[left,red] at (2.5,28) {$u+3v$};
    \node[right,red] at (97.5,44) {$2u+4v$};
    \node[right,red] at (97.5,17) {$u+v+1$};
    
    \node[below] at (50,0) {$f_1\colon [1, u + 3v] \to [u + v + 1, 2u + 4v]$};
\end{tikzoverlay}
\end{center}
  \caption{The triangles above show the faces and pairings that give rise to $f_1\colon [1, u + 3v] \to [u + v + 1, 2u + 4v]$. The arrows on the edges describe the orientation on their edges while the numbers on the pairings show their widths. The poles below are the edges aligned according to the integer intervals they correspond to. The strip in between the edges shows the pairing $f_1$. In particular, $f_1$ is orientation reversing.}
  \label{fig:fig2-1}
\end{figure}

Henceforth we will only consider a subset of $\{p'_i\}$, $\{f_1, g_1, h_1\}$ that consists of three special pairings. Let $f_1\colon [1, u + 3v] \to [u + v + 1, 2u + 4v]$, $g_1\colon [2v + 1, 3u + 3v] \to [u + 3v + 1, 4u + 4v]$, $h_1\colon [u + 3v + 1, 2u + 4v] \to [3u + 3v + 1, 4u + 4v]$ be pairings obtained as in Figures {\ref{fig:fig2-1}}, {\ref{fig:fig2-2}} and {\ref{fig:fig2-3}}. $f_1$ and $g_1$ are orientation reversing and $h_1$ is orientation preserving. Trim $f_1$ and $g_1$ to obtain $f_2\colon [1, u + 2v] \to [u + 2v + 1, 2u + 4v]$, $g_2\colon [2 v, 2u + 3v] \to [2u + 3v + 1, 4u + 4v]$. By composing $h_1$ with $g_2$ we can replace $h_1$ with $h_2\colon [u + 3v + 1, 2u + 4v] \to [2v + 1, u + 3v]$. $h_2$ is orientation reversing. Now $[1, 2v]$ and $[2u + 4v + 1, 4u + 4v]$ each lie in the domain and range of exactly one pairing so we can truncate $f_2$, $g_2$ to $f_3\colon [2v + 1, u + 2v] \to [u + 2v + 1, 2u + 2v]$, $g_3\colon [2u + 2v + 1, 2u + 3v] \to [2u + 3v + 1, 2u + 4v]$. Translating the domain and range of all pairings by $-2v$ we reach the interval $[1, 2u + 2v]$ with a collection of pairings $\{ f\colon [1, u] \to [u + 1, 2u], g\colon [2u + 1, 2u + v] \to [2u + v + 1, 2u + 2v], h\colon [1, u + v] \to [u + v + 1, 2u + 2v] \}$.

\begin{figure}[H]
\begin{center}
\begin{tikzoverlay}[height=100mm]{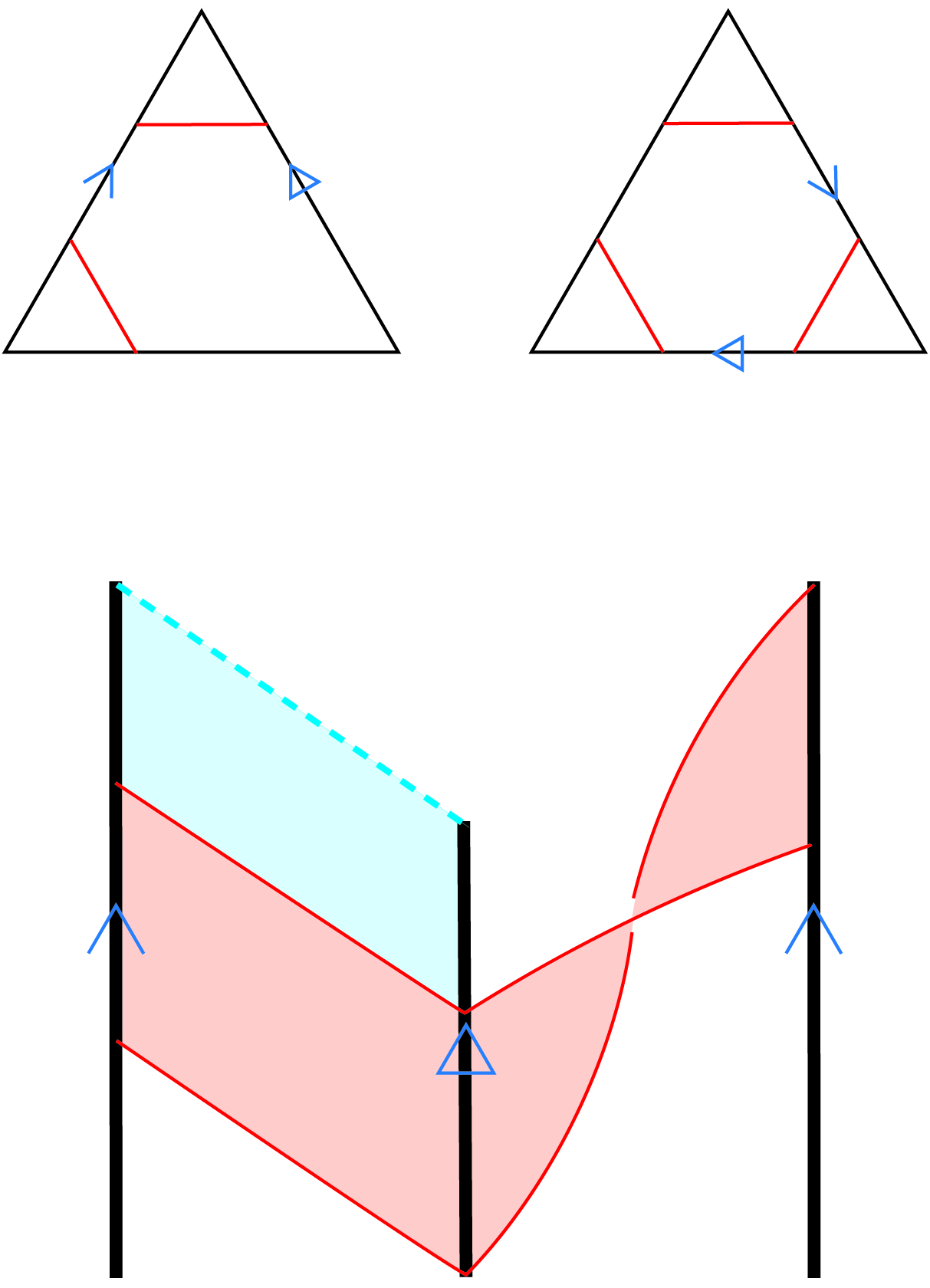}
\tikzstyle{every node}=[font=\small]
    \node[left] at (0,101) {$2$};
    \node[right] at (43,101) {$3$};
    \node[above] at (21.7,137) {$0$};
    \node[left] at (11,120) {e1};
    \node[right] at (33,120) {e3};
    \node[below] at (21.7,99.5) {F1 of tet0};
    \node[below,red] at (21.7,125) {$4u+2v$};
    \node[right,red] at (11,108) {$2v$};
    
    \node[left] at (57,101) {$1$};
    \node[right] at (100,101) {$3$};
    \node[above] at (78.4,137) {$0$};
    \node[right] at (89,120) {e1};
    \node[below] at (78.4,99.5) {e3};
    \node[below] at (78.4,95) {F2 of tet6};
    \node[below,red] at (78.4,125) {$u+3v$};
    \node[right,red] at (65,111) {$u+3v$};
    \node[left,red] at (89,106) {$3u+v$};
    
    \node[above] at (12.5,76) {e1};
    \node[above] at (87.5,76) {e1};
    \node[above] at (50,50) {e3};
    
    \node[left] at (12.5,2) {$1$};
    \node[left,red] at (12.5,27) {$2v+1$};
    \node[left,red] at (12.5,54) {$3u+3v$};
    \node[left] at (50,49) {$4u+2v$};
    \node[left,red] at (50,29) {$3u+v$};
    \node[right,red] at (87.5,75) {$4u+4v$};
    \node[right,red] at (87.5,48) {$u+3v+1$};
    
    \node[below] at (50,0) {$g_1\colon [2v + 1, 3u + 3v] \to [u + 3v + 1, 4u + 4v]$};
\end{tikzoverlay}
\end{center}
  \caption{Pairings that give rise to $g_1\colon [2v + 1, 3u + 3v] \to [u + 3v + 1, 4u + 4v]$. $g_1$ is orientation reversing.}
  \label{fig:fig2-2}
\end{figure}

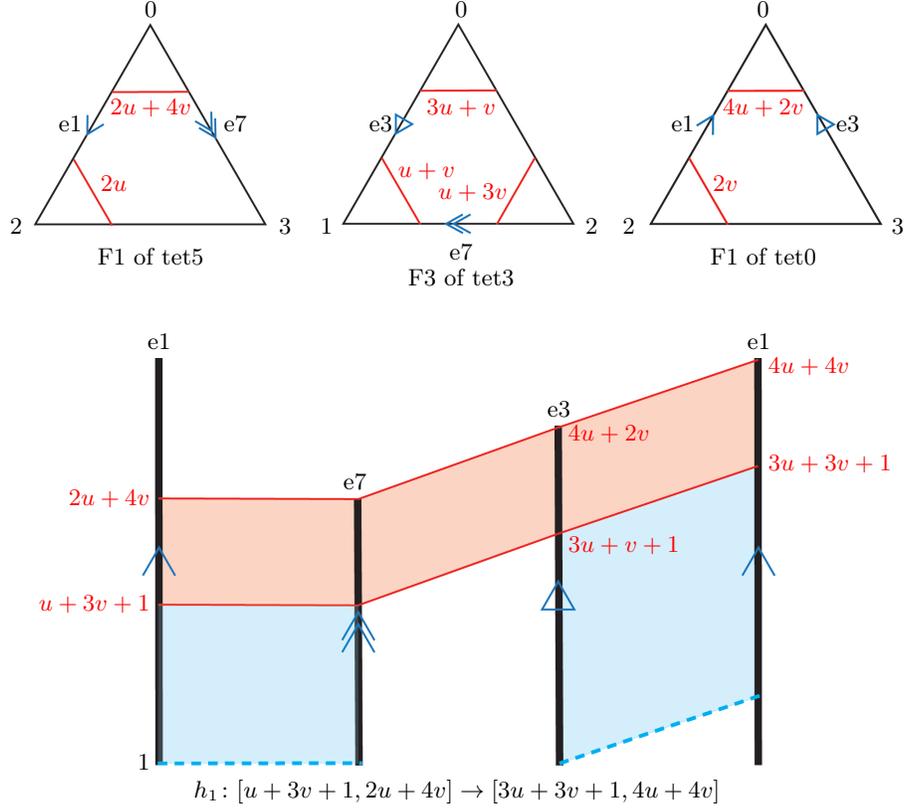
\begin{figure}[H]
\begin{center}
\begin{tikzoverlay}[height=100mm]{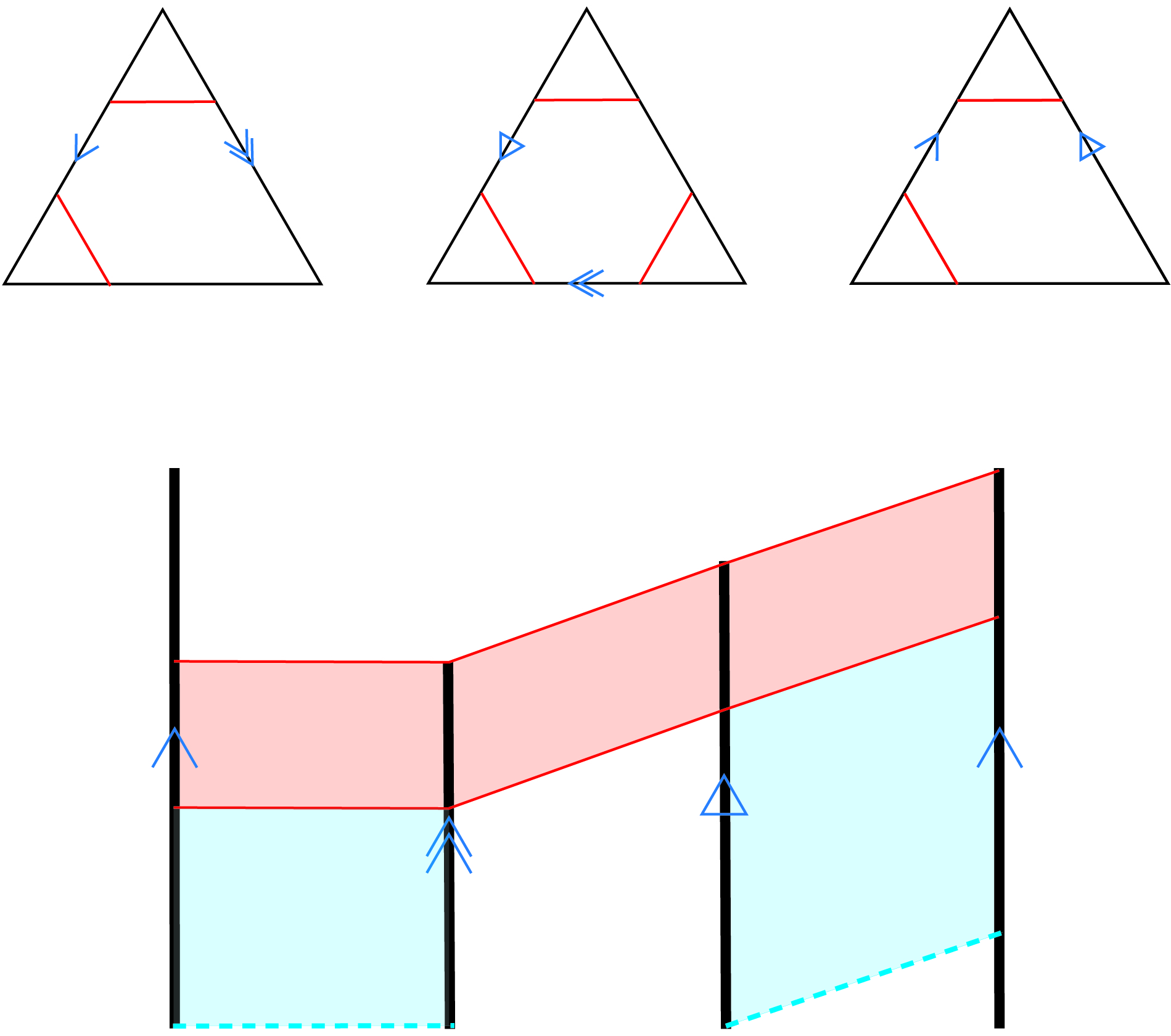}
\tikzstyle{every node}=[font=\small]
    \node[left] at (0,64) {$2$};
    \node[right] at (28,64) {$3$};
    \node[above] at (14,87.5) {$0$};
    \node[left] at (7,76) {e1};
    \node[right] at (21.5,76) {e7};
    \node[below] at (14,62.5) {F1 of tet5};
    \node[below,red] at (14,80) {$2u+4v$};
    \node[right,red] at (7,69) {$2u$};
    
    \begin{scope}[shift={(36,0)}]
    \node[left] at (36,64) {$2$};
    \node[right] at (64,64) {$3$};
    \node[above] at (50,87.5) {$0$};
    \node[left] at (43,76) {e1};
    \node[right] at (57.5,76) {e3};
    \node[below] at (50,62.5) {F1 of tet0};
    \node[below,red] at (50,80) {$4u+2v$};
    \node[right,red] at (43,69) {$2v$};
    \end{scope}

    \begin{scope}[shift={(-36,0)}]
    \node[left] at (72.5,64) {$1$};
    \node[right] at (100,64) {$2$};
    \node[above] at (86.5,87.5) {$0$};
    \node[left] at (79.5,76) {e3};
    \node[below] at (86.5,63) {e7};
    \node[below] at (86.5,60) {F3 of tet3};
    \node[below,red] at (86.5,80) {$3u+v$};
    \node[right,red] at (78,70.5) {$u+v$};
    \node[left,red] at (93,68) {$u+3v$};
    \end{scope}
    
    \node[above] at (15,48.5) {e1};
    \node[above] at (85.5,48.5) {e1};
    \node[above] at (38,32) {e7};
    \node[above] at (62,40.5) {e3};
    
    \node[left] at (15,1) {$1$};
    \node[left,red] at (15,19.5) {$u+3v+1$};
    \node[left,red] at (15,32) {$2u+4v$};
    \node[right,red] at (62,26.5) {$3u+v+1$};
    \node[right,red] at (62,39.5) {$4u+2v$};
    \node[right,red] at (85.5,47.5) {$4u+4v$};
    \node[right,red] at (85.5,36) {$3u+3v+1$};

    \node[below] at (50,0) {$h_1\colon [u + 3v + 1, 2u + 4v] \to [3u + 3v + 1, 4u + 4v]$};
\end{tikzoverlay}
\end{center}
  \caption{Pairings that give rise to $h_1\colon [u + 3v + 1, 2u + 4v] \to [3u + 3v + 1, 4u + 4v]$. Unlike $f_1$ and $g_1$, $h_1$ is orientation preserving.}
  \label{fig:fig2-3}
\end{figure}

\end{proof}

Now it remains to count the number of orbits.

\begin{claim} 
There are $gcd(u, v)$ orbits of $[1, 2u + 2v]$ under $\{f, g, h\}$.
\end{claim}

\begin{proof}

The three pairings $f, g, h$ are orientation reversing and $|f| = u, |g| = v, |h| = u + v$. $1$ is in the domain of $f$ whilst $2u + 2v$ is in the range of $g$ meaning $f$ lies on the left side of the interval and $g$ lies on the right side of the interval. Without loss of generality we assume $u \ge v$, if not we may simply reverse the order on the entire interval and apply the following proof. We will now describe a step by step orbit counting procedure with a sequence of pairings $\{ \{f_1, g_1, h_1 \}, \{f_2, g_2, h_2 \}, ..., \{f_N, g_N, h_N \} \}$ that starts with $\{f_1, g_1, h_1 \} = \{f, g, h \} $ and ends with a set of pairings $\{f_N, g_N, h_N \}$ where $|f_N| = |g_N|$. Geometrically, this process can be illustrated as in Figure {\ref{fig:fig3}}. The bands attached to the thickened interval in the middle illustrate the pairings.

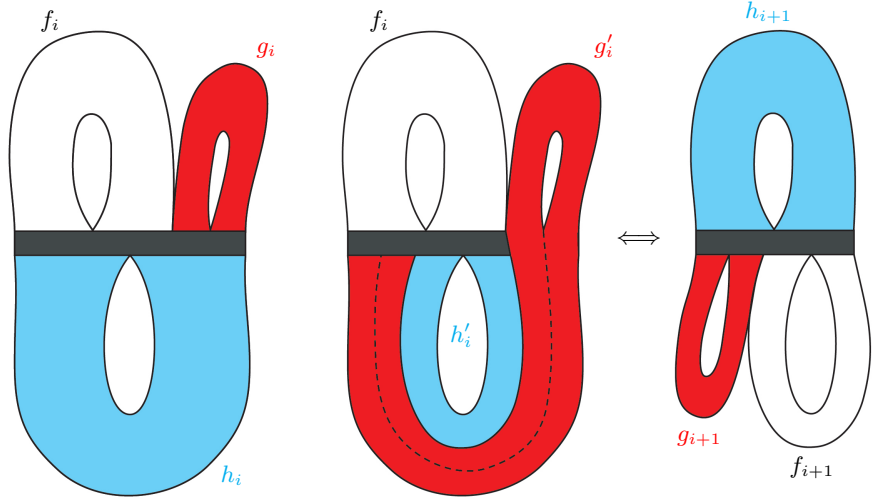
\begin{figure}[h]
\begin{center}
\begin{tikzoverlay}[width=115mm]{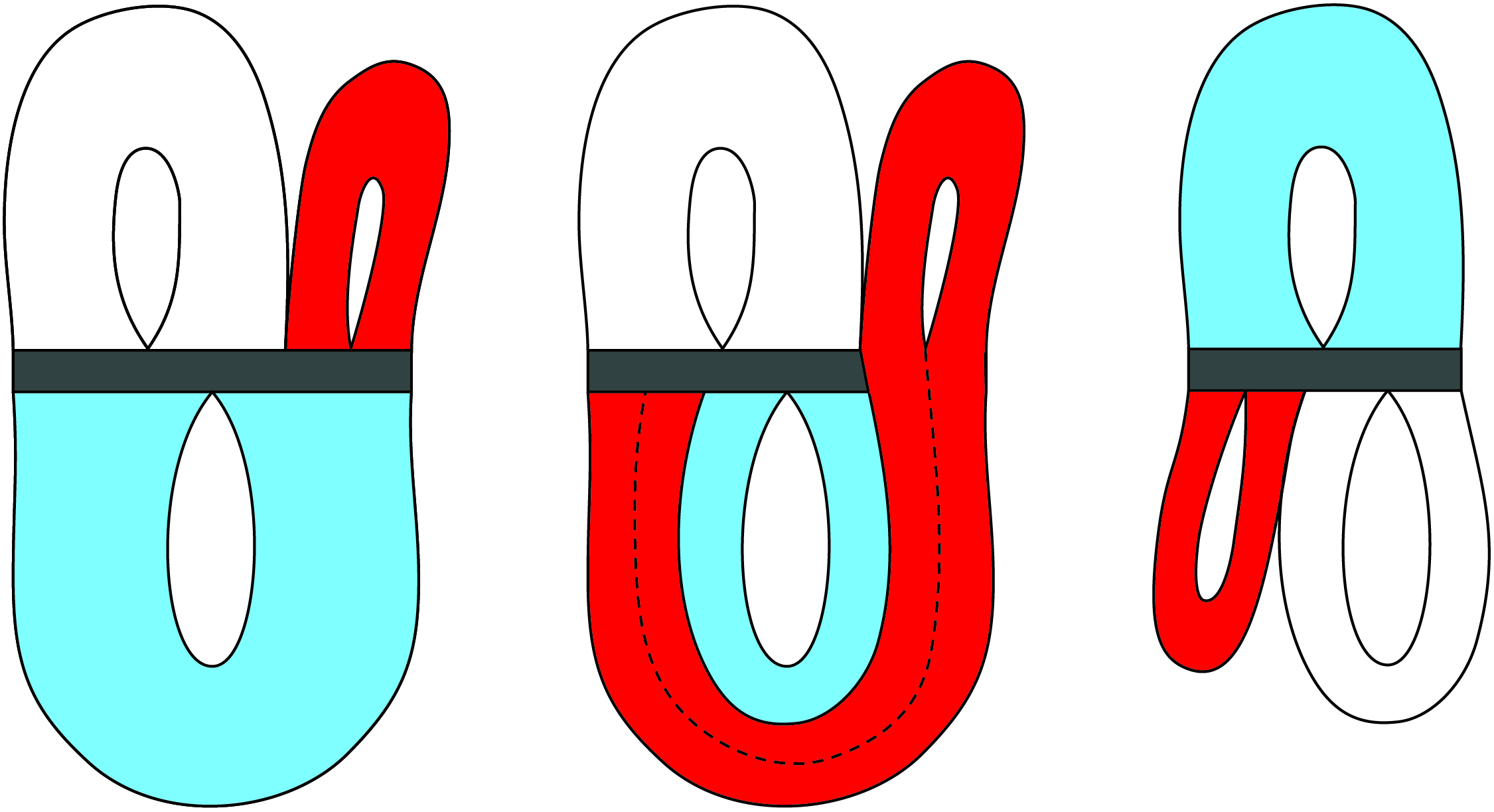}
\tikzstyle{every node}=[font=\small]
    \node[above] at (5,53) {$f_i$};
    \node[above,red] at (30,50) {$g_i$};
    \node[below,cyan] at (26,5) {$h_i$};
    
    \node[above] at (43,53) {$f_i$};
    \node[above,red] at (69,50) {$g_i'$};
    \node[cyan] at (52.5,19) {$h_i'$};
    
    \node[below] at (93,6) {$f_{i+1}$};
    \node[below,red] at (80,9) {$g_{i+1}$};
    \node[above,cyan] at (88,54) {$h_{i+1}$};
    
    \node at (73,30) {$\Longleftrightarrow$};
\end{tikzoverlay}
\end{center}
  \caption{Picture description of gcd calculation. The thick interval in the middle illustrates our interval while the bands attached are the pairings defined on subintervals.}
  \label{fig:fig3}
\end{figure}

\begin{enumerate}
  \item Transmit $g_i$ by $h_i$ to obtain $g_i'$.
  \item Truncate $h_i$ and peel off the domain and range of $g_i$ from $h_i$ to obtain $h_i'$.
  \item Assign $h_{i+1} = f_i$.
  \item Compare the widths of $g_i'$ and $h_i'$ then assign $f_{i + 1}$ with the pairing of larger width and $g_{i + 1}$ with the pairing of smaller width,
\end{enumerate}

In each step the pairings are of the form 
\begin{align*}
& f_i\colon \left[ a_i, \frac{a_i + b_i}{2} \right) \to \left( \frac{a_i + b_i}{2}, b_i \right], \\
& g_i\colon \left[ b_i, \frac{b_i + c_i}{2} \right) \to \left( \frac{b_i + c_i}{2}, c_i \right], \\ 
& h_i\colon \left[ a_i, \frac{a_i + c_i}{2} \right) \to \left( \frac{a_i + c_i}{2}, c_i \right], \quad a_i \le b_i \le c_i.
\end{align*}
The number of orbits are unchanged in this procedure and $|f_N| = |g_N|$ will give the number of orbits of $[1, 2u + 2v]$ under $\{f, g, h\}$ since  we have the pairings
\begin{align*}
& f_N\colon \left[ a_N, \frac{a_N + c_N}{4} \right) \to \left( \frac{a_N + c_N}{4}, \frac{a_N + c_N}{2} \right], \\
& g_N\colon \left[ \frac{a_N + c_N}{2}, \frac{3(a_N + c_N)}{4} \right) \to \left( \frac{3(a_N + c_N)}{4}, c_N \right] \\ 
& h_N\colon \left[ a_N, \frac{a_N + c_N}{2} \right) \to \left( \frac{a_N + c_N}{2}, c_N \right].
\end{align*}
Observe that in each step $|h_i| = |f_i| + |g_i|$, $|f_i| \ge |g_i|$. Moreover, $|h_{i + 1}| = |f_i|$ and $\{ |f_{i + 1}|, |g_{i + 1}| \} = \{ |g_i|, |f_i| - |g_i| \}$. The pairs of integers $(|f_i|, |g_i|)$ represent the sequence of integers that computes $gcd(|f_1|, |g_1|) = gcd(u, v)$ by Euclid's algorithm. Hence the number of orbits $|f_N| = |g_N| = gcd(u, v)$.
  
\end{proof}

By assumption $gcd(u, v) = 1$ so we are done by Claim $1$.

\end{document}